\documentclass[final]{siamltex}
\usepackage{epsfig}
\usepackage{amssymb}
\usepackage{amsmath, amssymb}
\usepackage{graphics,color}

\newtheorem{example}[theorem]{Example}

\def\qed{\hfill $\Box$\medskip}
\def\IC{{\mathbb C}}
\def\IR{{\mathbb R}}

\def\cR{{\mathcal R}}

\def\({\left (}
\def\){\right )}

\def\conv{{\rm conv}\,}

\begin{document}
\openup .8\jot

%%%%%%%%%%%%%%%%%%%%%%%%%%%%%%%%%%%%%%%%%%%%%%%%%%%%%%%%%%%%%%%%%%%%%%%%%%%%%%%%%%%%%%%%%%%%%%%%%%%%%%%%%%%%%%%%%%%%%%%%%%%%%%%%%%%%%%%%%%%%%%%%%%%%%%%%%%%%%%%%%%%%%%%%%%%%%%%%%%%%%%%%%%%%%%%%%%%%%%%%%%%%%%%%%%%%%%%%%%%%%%%%%%%%%%%%%%%%%%%%%%%%%%%%%%%%%%%%%%%%%%%%%%%%%%%%%%%%%%%%%%%%%%%%%%%%%%%%%%%%%%%%%%%%%%%%%%%%%%%%%%%%%%%%%%%%%%%%%%%%%%%%%%%%%%%%%%%%%%%%%%%

\title{
Determinantal and eigenvalue inequalities\\
for matrices with 
numerical ranges in a sector}

\author{
Chi-Kwong Li\thanks{
Department of Mathematics, College of William and Mary, Williamsburg, VA 23187, USA
(Email: ckli@math.wm.edu).
Li is an honorary professor of the
Taiyuan University of Technology (100 talents program),
the University of Hong Kong, and also the
Shanghai University. His research was partially supported
by a USA NSF grant DMS 1007768 and a HK RGC grant.
}
\and
Nung-Sing Sze\thanks{
Department of Applied Mathematics,
The Hong Kong Polytechnic University,
Hung Hom, Hong Kong
(Email: raymond.sze@polyu.edu.hk).
Research of Sze was partially supported by
a HK RGC grant PolyU 502512.
}
}

\maketitle

\begin{abstract}
Let $A = \begin{pmatrix}
A_{11} & A_{12} \cr A_{21} & A_{22}\cr\end{pmatrix} \in M_n$, where
$A_{11} \in M_m$ with $m \le n/2$,
be such that the numerical range of $A$ lies in the set
$\{e^{i\varphi} z \in \IC: |\Im z| \le (\Re z) \tan \alpha\}$,
for some $\varphi \in [0, 2\pi)$ and 
$\alpha \in [0, \pi/2)$.
We obtain the optimal containment region for the 
generalized eigenvalue $\lambda$ satisfying 
$$\lambda \begin{pmatrix}
A_{11} & 0 \cr 0 & A_{22}\cr\end{pmatrix}  x = 
\begin{pmatrix}0 & A_{12} \cr A_{21} & 0\cr\end{pmatrix} x \quad 
\hbox{ for some nonzero } x \in \IC^n,$$
and  the optimal eigenvalue containment region of
the matrix $I_m - A_{11}^{-1}A_{12} A_{22}^{-1}A_{21}$
in case $A_{11}$ and $A_{22}$ are invertible.
From this result, one can show
$|\det(A)| \le \sec^{2m}(\alpha) |\det(A_{11})\det(A_{22})|$.
In particular,  if $A$ is a accretive-dissipative matrix, then
$|\det(A)| \le 2^m |\det(A_{11})\det(A_{22})|$.
These affirm some conjectures of Drury and Lin.
\end{abstract}

\begin{AMS}
15A45.
\end{AMS}

\begin{keywords}
Numerical ranges,  
eigenvalues, determinantal inequality, 
accretive-dissipative matrix.
\end{keywords}
\section{Introduction}

Let $M_n$ be the set of $n\times n$ complex matrices.
Suppose
\begin{equation}\label{formA}
A = \begin{pmatrix}
A_{11} & A_{12} \cr A_{21} & A_{22}\cr\end{pmatrix} \in M_n \quad 
\hbox{ with } A_{11} \in M_m, \quad m \le n/2.
\end{equation}
In connection to the study of the growth factor in
Gaussian elimination, researchers considered optimal 
(smallest) $\gamma > 0$ such that 
$$|\det(A)| \le \gamma |\det(A_{11})\det(A_{22})|;$$
see \cite{D,GI,GI2,H,I,Lin} and their references.
The well-known Fischer inequality asserts that
$$\det(A) \le \det(A_{11})\det(A_{22}) \quad \hbox{ if } A \hbox{ is positive 
semi-definite}.$$
In \cite{I}, it was shown that if
$A$ is accretive and dissipative, i.e.,   $A+A^*$ and $i(A^*-A)$ are
positive semi-definite, then  
$$|\det(A)| \le \gamma |\det(A_{11})\det(A_{22})| 
\qquad \hbox{ with } \gamma = 3^m;$$
in \cite{Lin}, the bound was improved to
$$\gamma = \begin{cases}2^{3m/2} & \hbox{ if } m\le n/3,\\[1mm]
2^{n/2} & \hbox{ if } n/3 < m \le n/2.\cr\end{cases}$$
The author in \cite{Lin} further proposed the following.

\medskip\noindent
{\bf Conjecture 1} \it Suppose $A$ is accretive-dissipative. Then 
\begin{equation}
|\det(A)| \le 2^m |\det(A_{11})\det(A_{22})|.
\end{equation}
\rm

\medskip
The numerical range of a matrix $L \in M_n$ is defined by
$$W(L) = \{x^*Lx: x \in \IC^n, \ x^*x = 1\}.$$
For any $\alpha \in [0,\pi/2)$, let 
\begin{equation}\label{formSa}
S_{\alpha} = \{z \in \IC: |\Im z| \le (\Re z) \tan \alpha\}.
\end{equation}
A subset of $\IC$ is a sector of half angle $\alpha$ if
it is of the form $\{e^{i\varphi}z: z \in S_\alpha\}$
for some $\varphi \in [0, 2\pi)$.

In \cite{D}, the author proved that if $W(A)$ is a subset of 
a sector of half angle $\alpha \in [0, \pi/(2m))$,   then   
$$|\det(A)| \le \sec^2(m\alpha) |\det(A_{11})\det(A_{22})|.$$
He further proposed the following.

\medskip\noindent
{\bf Conjecture 2} \it
If $W(A)$ is a subset of 
a sector of half angle $\alpha \in [0, \pi/2)$, then
\begin{equation} \label{conj1}
|\det(A)| \le \sec^{2m}(\alpha) |\det(A_{11})\det(A_{22})|.
\end{equation}
Moreover, if $A_{11} \in M_1$ is nonzero and $A_{22}\in M_{n-1}$ 
is invertible, then
\begin{equation} \label{conj02}
\det(A)/(\det(A_{11})\det(A_{22})) \in\{r e^{i2\phi}: 
0 \le r \le 2(\cos (2\phi) - \cos (2\alpha))/\sin^2(2\alpha),
-\alpha \le \phi \le \alpha\}.
\end{equation}
\rm

\medskip
We will affirm Conjectures 1 and 2 via the study of
the following generalized eigenvalue problem
\begin{equation}\label{geneig}
\lambda \begin{pmatrix}A_{11} & 0 \cr 0 & A_{22} \cr\end{pmatrix} x = 
\begin{pmatrix}0& A_{12} \cr A_{21} & 0 \cr\end{pmatrix} x
\quad \hbox{ for some nonzero } x \in \IC^n. 
\end{equation}
In Section 2, we prove Theorem \ref{main}
providing the optimal eigenvalue containment region for
those $\lambda$ satisfying (\ref{geneig}),
and  the optimal eigenvalue containment region of
the matrix 
$$A_{11}^{-1}A_{12} A_{22}^{-1}A_{21}$$
in case $A_{11}$ and $A_{22}$ are invertible.
Using the theorem, one can readily verify Conjectures 1 and 2;
see Corollaries \ref{cor1} and \ref{cor2}.

\section{Results, proofs, and remarks}

In this section, we will always assume that 
$A$ has the form in (\ref{formA}),
and refer to a subset of $\IC$ as
a sector of half angle $\alpha$ if
it is of the form $\{e^{i\varphi}z: z \in S_\alpha\}$
for some $\varphi \in [0, 2\pi)$ and $S_\alpha$ defined in (\ref{formSa}).

We begin with several lemmas.
 
\begin{lemma} \label{lemf} For any $\phi \in [-\pi/2,\pi/2]$,
the function  $f: (0, \pi/2) \rightarrow \IR$
defined by 
$$f(\theta) =  (\cos (2\phi) - \cos (2\theta)) /\sin^2 (2\theta)$$
is increasing.
\end{lemma}

\it Proof. \rm By direct verification. \qed

\begin{lemma} \label{lem}
Suppose $z_1 = r_1e^{i\theta_1}, z_2 = r_2e^{i\theta_2} \in \IC$ with
$r_1, r_2 \in (0,\infty)$, $\pi/2 > \theta_1 \ge \theta_2 > -\pi/2$,
and $(z_1+z_2)/2 = e^{i\psi}$.
Let $2\theta = \theta_1-\theta_2$ and $2\phi = \theta_1+\theta_2-2\psi$.
Then 
$$r_1r_2 = 2(\cos(2\phi)-\cos(2\theta))/\sin^2(2\theta).$$
\end{lemma}

\it Proof. \rm
Consider the triangle $T$ with vertices $0, z_1, z_2$. 
Because $e^{i\psi}$ is the midpoint of the side
joining the vertices $z_1$ and $z_2$,
 $T$ can be divided into two triangles $T_1$
and $T_2$ with equal areas, where $T_1$ has 
vertices $0, e^{i\psi}, r_1e^{i\theta_1}$, and $T_2$ has vertices
$0, e^{i\psi}, r_2e^{i\theta_2}$.
Thus, 
$$r_1r_2\sin(2\theta) = r_1 \sin(\theta_1-\psi) + r_2\sin(\psi-\theta_2 )
\quad\hbox{ and } \quad 
r_1 \sin(\theta_1-\psi) = r_2 \sin(\psi-\theta_2).$$
It follows that
\begin{eqnarray*}
(r_1r_2\sin(2\theta))^2 
&=& (r_1\sin(\theta_1-\psi) + r_2\sin(\psi-\theta_2))^2\cr
&=& 4 r_1r_2\sin(\theta_1-\psi)\sin(\psi-\theta_2)\cr
&=& 2 r_1r_2(\cos (2\phi) - \cos(2\theta)).
\end{eqnarray*}
Hence, $r_1r_2 = 2(\cos (2\phi) -  \cos(2\theta))/\sin^2(2\theta)$.
\qed

We will also use some basic facts about
the numerical range; for example, see \cite[Chapter 1]{HJ}.

\begin{lemma} \label{wa} Let $L \in M_n$.
\begin{enumerate}
\item If $U \in M_n$ is unitary, then $W(L) = W(U^*LU)$.
\item The set $W(L)$ is compact and convex.
\item If $\tilde L$ is a principal submatrix of $L$, then 
$W(\tilde L) \subseteq W(L)$.
\item If $L = L_1 \oplus L_2$, then
$W(L) = \conv(W(L_1)\cup W(L_2))$, where $\conv(S)$ denotes the convex hull
of the set $S$. 
\item If $L$ is normal, then $W(L)$ is the
convex hull of its eigenvalues.
\item If $x \in \IC^n$ is a unit vector such that
$\mu = x^*Lx$ is a boundary point of 
$W(L)$ with more than one support line,
then $Lx = \mu x$ and $L^*x = \overline{\mu} x$.
\end{enumerate}
\end{lemma}

\medskip
\begin{theorem} \label{main}
Let 
$A = \begin{pmatrix}A_{11} & A_{12} \cr A_{21} & A_{22}\cr
\end{pmatrix} \in M_n$
with $A_{11} \in M_m$ be such that $m \le n/2$, and 
$W(A)$ be a subset of a sector of half angle $\alpha \in [0, \pi/2)$.

\begin{itemize}
\item[{\rm (a)}] Suppose $A_{11} \oplus A_{22}$ is singular, and 
$x \in \IC^n$ is a nonzero vector in its kernel. Then 
$Ax = 0$ and {\rm (\ref{geneig})} holds for every $\lambda \in \IC$
with this nonzero vector $x$.

\item[{\rm (b)}] Suppose $A_{11}$ and $A_{22}$ are invertible.
If $\lambda\in \IC$ satisfies {\rm (\ref{geneig})}, then
$$\lambda^2 \in \cR = \begin{cases} 
\left\{1-re^{i2\phi}: 0 \le r \le \frac{2 (\cos (2\phi) - \cos (2\alpha))}
{\sin^2 (2\alpha)}, -\alpha \le \phi \le \alpha \right\} & \hbox{ if } \alpha > 0,\cr
\,[0,1] & \hbox{ if } \alpha = 0. \cr\end{cases}$$
\iffalse
where we use the convention that
$ \frac{2 (\cos (2\phi) - \cos (2\alpha))}{\sin^2 (2\alpha)} = 1$
if $\alpha = 0$. \fi
Moreover, for every eigenvalue
$\mu$ of the matrix
$A_{11}^{-1} A_{12} A_{22}^{-1} A_{21}$, 
there is $\lambda \in \IC$ satisfying {\rm (\ref{geneig})} so that
$\mu = \lambda^2$ lies in the region $\cR$.
\end{itemize}
\end{theorem}
 
\it Proof. \rm 
Without loss of generality, we may assume that $W(A) \subseteq S_\alpha$.
Let 
$$B_1 = \begin{pmatrix}A_{11} & 0 \cr 0 & A_{22}\cr\end{pmatrix} 
\qquad \hbox{ and } \qquad
B_2 = \begin{pmatrix}0 & A_{12} \cr A_{21} & 0 \cr\end{pmatrix}.$$ 

(a) Suppose $B_1$ is singular and $x \in \IC^n$ is nonzero
such that $B_1 x = 0$. By Lemma \ref{wa}, for $y = x/\|x\|$,
$$0 = y^*B_1 y \in \conv(W(A_{11}\oplus A_{22})) \subseteq W(A)$$
so that 0 is a boundary point of $W(A)$ with more than one support line.
Thus, $Ax = 0$, and $\lambda B_1 x = 0 = (A-B_1)x = B_2 x$
for any $\lambda\in \IC$.

(b) Assume $A_{11}$ and $A_{22}$ are invertible.
Suppose $\lambda B_1 x = B_2 x$ for some nonzero unit vector $x\in \IC^n$.
Let $\xi_1 = x^*B_1x$ and $\xi_2 = x^*B_2x$.
Then $\xi_1 \in W(B_1)$ and $\xi_2 \in W(B_2)$.
We see that $\xi_1 + \xi_2 \in W(B_1+B_2) = W(A)$ and
$\lambda \xi_1 = \xi_2$.
Note that $A = B_1 + B_2 = Q^*(B_1 - B_2) Q$
with $Q = I_m \oplus (-I_{n-m})$, and hence $W(A) = W(B_1-B_2)$
contains $x^*(B_1-B_2)x = \xi_1 - \xi_2$. So,
$\xi_1 \pm \xi_2 \in W(A)$, which is  a subset of 
$S_{\alpha}$ by our assumption.

Observe that $\xi_1 \ne 0$. Otherwise,  by Lemma \ref{wa}
$$0 \in W(B_1) = \conv(W(A_{11}) \cup W(A_{22})) \subseteq W(A)$$
so that 0 is a boundary point of $W(B_1)$ 
with more than one support line implying that $B_1$ is singular,
which contradicts our assumption.

Without loss of generality, assume that 
$\Im(\lambda) = \Im(1+\xi_2/\xi_1) \ge 0$. 
Let
$$z_\pm = r_\pm e^{i\theta_\pm} = 1\pm \xi_2/\xi_1
= 1 \pm \lambda \quad\hbox{ with} \quad r_\pm \ge 0,
\quad \theta = \frac{1}{2} (\theta_+ - \theta_-)
\quad\hbox{and}
\quad \phi = \frac{1}{2} (\theta_+ + \theta_-).$$
If $\xi_1 = |\xi_1| e^{i\omega}$, 
then $(\xi_1 \pm \xi_2) = \xi_1(1\pm \xi_2/\xi_1)$ 
has arguments $\theta_\pm + \omega \in [-\alpha,\alpha]$
as $\xi_1\pm \xi_2 \in W(A) \subseteq S_\alpha$.
Note also that $(z_1+z_2)/2 = 1$ has argument 0. It follows that 
$-\alpha - \omega \le \theta_- \le 0 \le \theta_+ \le \alpha - \omega$.
So $0 \le \theta \le \alpha$.
\iffalse and $-\alpha \le \psi \le \alpha$,
where $\omega$ is the argument of $\xi_1$ with $-\alpha\le \omega \le \alpha$.\fi
Applying Lemma \ref{lem} with $(r_1e^{i\theta_1},r_2e^{i\theta_2}) 
= (r_+e^{i\theta_+},r_-e^{i\theta_-})$
and $\psi = 0$,
%and $(\psi,\phi) = (0,(\theta_+ + \theta_-)/2)$, 
we have 
$$
r_+r_- = \frac{2 ( \cos (2\phi) - \cos (2 \theta))}{\sin^2 (2\theta)}
\le \frac{2 (\cos (2\phi) - \cos (2 \alpha))}{\sin^2 (2\alpha)},
$$
where the  inequality follows from Lemma \ref{lemf}.
As a result, 
$$1-\lambda^2 = z_+ z_- 
= r_+ r_- e^{i(\theta_++\theta_-)} = r_+r_-e^{i2\phi}$$
lies in the region
$$\tilde R = \left\{r e^{i2\phi}: 0 \le r \le 
\frac{2 ( \cos (2\phi) - \cos (2\alpha))}{\sin^2 (2\alpha)}, 
-\alpha \le \phi \le \alpha \right\}.$$

Suppose $B = B_1^{-1}B_2$.
Then
$$B^2 = \begin{pmatrix} A_{11}^{-1}A_{12} A_{22}^{-1}A_{21} & 0 \cr
0 & A_{22}^{-1}A_{21}A_{11}^{-1}A_{12} \cr\end{pmatrix}.$$
If $A_{11}^{-1}A_{12} A_{22}^{-1}A_{21}$ has eigenvalues
$\mu_1, \dots, \mu_m$, then  $A_{22}^{-1}A_{21}A_{11}^{-1}A_{12}$
has eigenvalues $\mu_1, \dots, \mu_m$ together with $n-m$ zeros.
So, we may assume that $B$ has eigenvalues $\lambda_1, \dots, \lambda_n$
such that $\lambda_j^2 = \lambda_{m+j}^2 = \mu_j$ for $j = 1, \dots, m$
and $\lambda_\ell = 0$ for $\ell = 2m+1, \dots, n$.
Note that $\lambda$ is an eigenvalue of $B$ if and only if
$\lambda$ satisfies (\ref{geneig}) for some nonzero $x \in \IC^n$.
The second assertion of (b) follows.
\qed

The containment region in Theorem \ref{main}\,(b) is optimal as shown in 
the following.

\begin{example} \label{ex} \rm
Let $\lambda \in \IC$ be such that $\lambda^2 \in \cR$ in Theorem \ref{main}, 
i.e., $1-\lambda^2 = re^{i2\phi}$ with
$r \in \left[0, {2(\cos (2\phi) - \cos (2\alpha))}/{\sin^2(2\alpha)}\right]$
for some $\alpha \in [0, \pi/2)$.

Suppose $r > 0$.
By Lemma \ref{lemf}, there is 
$\theta \in (|\phi|,\alpha]$ satisfying
$r = 2(\cos (2\phi) - \cos (2\theta))/\sin^2(2\theta)$, 
here we set
$2(\cos (2\phi) - \cos (2\theta))/\sin^2(2\theta) = 1$ if $\theta=0$ 
which will imply $\phi=0$.
Let
$$A =  \left(
I_m \otimes \begin{pmatrix}e^{-i\phi}\ & a+ib \cr 
a+ib & e^{-i\phi} \cr\end{pmatrix} \right) \oplus (e^{-i\phi} I_{n-2m})
$$
with $a = -\cot\theta\sin\phi$ and $b = \tan\theta\cos\phi$ so that 
$$|a| \le |\cot\phi\sin\phi| = \cos\phi,  \
\quad b \ge |\tan\phi\cos\phi| = |\sin\phi|,$$
\begin{equation} \label{ab}
\frac{-\sin\phi+b}{\cos \phi+a} = \tan\theta, 
\quad \hbox{ and } \quad
\frac{-\sin\phi-b}{\cos \phi - a} = -\tan\theta.
\end{equation}
Then $A$ is normal, and by (\ref{ab}) the eigenvalues 
of $A$ has the form
$e^{-i\phi} + (a+ib) = r_1 e^{i\theta}$, $r_1 \ge 0$, with multiplicity $m$,
$e^{-i\phi} - (a+ib) = r_2 e^{-i\theta}$, $r_2 \ge 0$, with multiplicity $m$, 
and $e^{-i\phi}$ with multiplicity $n-2m$, all in $S_\alpha$.
By Lemma \ref{wa}, $W(A) \subseteq S_\alpha$.
Moreover,
$\lambda = \pm (a+ib)e^{i\phi}$ satisfy (\ref{geneig}),
and 
$$1-\lambda^2 = \det\begin{pmatrix}1\ & (a+ib)e^{i\phi} \cr 
(a+ib)e^{i\phi}& 1 \cr\end{pmatrix}
= e^{2i\phi}\det\begin{pmatrix}e^{-i\phi}\ & a+ib \cr 
a+ib & e^{-i\phi} \cr\end{pmatrix} 
= e^{2i\phi} r_1 e^{i\theta} r_2 e^{-i\theta}
= r_1r_2e^{2i\phi}.$$
Applying Lemma \ref{lem}
to $r_1 e^{i\theta}, r_2 e^{-i\theta}$
so that $\theta_1 = -\theta_2 = \theta$ and $\psi = -\phi$, and using the fact that 
$e^{-i\phi}$ is the midpoint of the line segment
joining $r_1 e^{i\theta}, r_2 e^{-i\theta}$, we have
$$r_1r_2 = 2(\cos (2\phi) -  \cos(2\theta))/\sin^2(2\theta) = r.$$

Suppose $r = 0$. One can verify directly that
the matrix $$A =  \left(
I_m \otimes \begin{pmatrix} 1 & 1 \cr 
1 & 1 \cr\end{pmatrix} \right) \oplus I_{n-2m}
$$
satisfies $W(A) = [0,2] \subseteq S_\alpha$ and has
generalized eigenvalues $\lambda \in \{1,-1\}$ so that $1 - \lambda^2 = 0 = r$.
\qed
\end{example}

Using the notation of Theorem \ref{main}, 
we see that if $A_{11}$ and $A_{22}$ are invertible, then 
all eigenvalues of the matrix 
$C = I_m - A_{11}^{-1} A_{12} A_{22}^{-1} A_{21}$ lie in the set 
$\left\{re^{i2\phi}: 0 \le r \le \frac{2 (\cos (2\phi) - \cos (2\alpha))}
{\sin^2 (2\alpha)}, |\phi| \le  \alpha \right\}.$
Thus, the spectral radius of $C$ is bounded by
$$\max_{|\phi|\le \alpha} \left\{\frac{2 (\cos (2\phi) - \cos (2\alpha))}
{\sin^2 (2\alpha)} \right\}= \frac{2 (1 - \cos (2\alpha))}
{\sin^2 (2\alpha)} = \sec^2(\alpha),$$
and hence
$$|\det(A)| = |\det(A_{11})\det(A_{22})\det(C)| 
\le \sec^{2m}(\alpha) |\det(A_{11})\det(A_{22})|.$$ 
By continuity, one can remove the 
invertibility assumption on 
$A_{11}$ and $A_{22}$. We have the following corollary 
affirming Conjecture 2.

\begin{corollary} \label{cor1}
Let 
$A = \begin{pmatrix}A_{11} & A_{12} \cr A_{21} & A_{22}\cr\end{pmatrix} \in M_n$
with $A_{11} \in M_m$ such that $m \le n/2$, and 
$W(A)$ be a subset of a sector of half angle $\alpha \in [0, \pi/2)$.
Then 
$$|\det(A)| 
\le \sec^{2m}(\alpha) |\det(A_{11})\det(A_{22})|.$$ 
If $A_{11}$ and $A_{22}$ are invertible, then the eigenvalues of the matrix
$C = I_m - A_{11}^{-1} A_{12} A_{22}^{-1} A_{21} \in M_m$ lies in the region 
%$$\tilde \cR = \{re^{i2\phi}: 
%0 \le r \le \sec^2(\alpha-|\phi|), |\phi| \le \alpha\}.$$
$$\tilde \cR = 
\left\{re^{i2\phi}: 0 \le r \le \frac{2 (\cos (2\phi) - \cos (2\alpha))}
{\sin^2 (2\alpha)}, -\alpha \le \phi \le  \alpha \right\}.$$
\end{corollary}

Suppose $A$ is accretive-dissipative. Then 
$W(e^{-i\pi/4} A) \subseteq S_{\pi/4}$. Applying Corollary \ref{cor1}
with $\alpha = \pi/4$, we have the following result
affirming Conjecture 1 and verifying a comment in \cite{D} (after Conjecture 0.2).

\begin{corollary} \label{cor2}
Let 
$A = \begin{pmatrix}A_{11} & A_{12} \cr A_{21} & A_{22}\cr\end{pmatrix} \in M_n$ be 
accretive-dissipative with $A_{11} \in M_m$ such that $m \le n/2$. Then
$$|\det(A)| \le 2^m |\det(A_{11})\det(A_{22})|.$$ 
If $A_{11}$ and $A_{22}$ are invertible, then the eigenvalues of the matrix
$C = I_m - A_{11}^{-1} A_{12} A_{22}^{-1} A_{21} \in M_m$ lies in the set
$$\{ z\in \IC: |z-1| \le 1\}.$$
\end{corollary}

The bound in Corollary \ref{cor1} is best possible by 
Example \ref{ex} with $(\phi,\theta) = (0,\alpha)$.
In such a case, $(e^{-i\phi},a+bi) = (1,i\tan \alpha)$, 
$\det(A) = (1+\tan^2\alpha)^{m} = \sec^{2m}(\alpha)$ and 
$1 = \det(A_{11}) = \det(A_{22})$.
Furthermore, letting $\alpha= \pi/4$, 
we see that the bound in Corollary \ref{cor2} is also best possible.

\medskip\noindent
{\bf Acknowledgment}

We thank Dr.\ Minghua Lin for drawing our attention to this 
interesting topic,
and Professor S.W. Drury for his preprint \cite{D}.
Upon the completion of our note, we learned that
Professor Drury had verified conjecture 2 for the case
when $\alpha \in [0, \pi/4]$ by a different method in \cite{D2}.

Thanks are also due to Dr.\ Zejun Huang for some helpful discussion.

\end{document}